\documentclass[11pt]{article}
\usepackage{amsmath,amssymb}

\def\epsilon{\varepsilon}

\def\phi{\varphi}

\newtheorem{theorem}{Theorem}[section]
\newtheorem{lemma}[theorem]{Lemma}

\newtheorem{remark}[theorem]{Remark}

\def\ra{A^{\!\scriptscriptstyle1\!/\kern-0.8pt 2}}

\def\R{{\mathbb R}}

\newenvironment{Proof}{\removelastskip\par\medskip
\noindent{\em Proof.} \rm}{\penalty-20\null$\square$\par\medbreak}

\newenvironment{Proofy}{\removelastskip\par\medskip
\noindent{\em Proof of Theorem~\ref{th:intin}.} 
\rm}{\penalty-20\null$\square$\par\medbreak}

\title{\bf Energy decay   for evolution  equations 
\\ with glassy type memory}

\author{Paola Loreti
\thanks{Dipartimento di Scienze di Base e Applicate per l'Ingegneria,
Sapienza Universit\`a di Roma,
Via Antonio Scarpa 16, 00161 Roma (Italy); e-mail: 
$<$paola.loreti@uniroma1.it$>$ }
\and Daniela Sforza
\thanks{Dipartimento di Scienze di Base e Applicate per l'Ingegneria, 
Sapienza Universit\`a di Roma,
Via Antonio Scarpa 16, 00161 Roma (Italy); e-mail: 
$<$daniela.sforza@uniroma1.it$>$ }}

\begin{document}
\date{}

\maketitle

\begin{abstract}
In this paper, we address the question of estimating the energy decay of integro-differential evolution equations with glassy memory.
This class of memory kernel was not analyzed in previous studies. Moreover, a detailed analysis provides an explicit estimate of the connection between the kernel function's decay constant and the energy's decay constant.   
\end{abstract}


\bigskip
\noindent

\section{Introduction}
In this paper, we establish decay estimates of the energy for general second-order integro-differential evolution equations. These estimates depend on the properties of the convolution kernels, which do not satisfy one of the main assumptions in \cite{ACS}. Nevertheless, the solution's energy still decays exponentially.
Moreover, a detailed analysis provides an explicit estimate of the connection between the kernel function's decay constant and the energy's decay constant. 
We use a method based on multiplier techniques.

As an example, we can consider the viscoelastic wave equation with glassy memory and a Petrovsky-type system. In this introduction, we will focus on the wave equation with memory to illustrate the main steps.
However, the results will be developed within an abstract and general framework in later sections.

To contextualize our findings, we will briefly review some basic facts.

Let $\Omega$ be a bounded, open and connected  set of $\R^N$ with boundary 
$\Gamma$ of class $C^2$. %Let  $T>0$
Let us consider the homogeneous wave equation, which models vibrations in elastic strings or plates:
\begin{equation}\label{eq:ondaI}
\begin{cases}
u_{tt} -\triangle u=0  \quad {\text {in}}\ \R\times\Omega, 
\\
u(0,x)=u_0(x)\quad u_t(0,x)=u_1(x)\quad {\text {in}}\ \Omega, 
\\
u(t,x)=0 \quad {\text {on}}\ \R\times\Gamma.
\end{cases}
\end{equation}
For $u_0\in H^1_0(\Omega)$ and $u_1\in L^{2}(\Omega)$
\eqref{eq:ondaI} has a unique weak solution $u$ belonging to 
$
C(\R ;H_0^1(\Omega))\cap C^1(\R ;L^{2}(\Omega)).
$
Accordingly, as well known, the energy of a weak solution can be defined as
\begin{equation*}
E(t)=\frac{1}{2}\int _\Omega u_t^2+|\nabla u|^2\, dx
\end{equation*}
and it is conserved, because, multiplying the equation by $u_t$ and taking into account the boundary condition, one obtains
$E'(t)=0$ for all $t$, that is the system is conservative.

As is well known, the viscoelastic behavior of materials combines elastic and viscous components, often modeled as linear combinations of springs and dashpots. These materials are widely applied in fields like medicine, industry, and biology, and are particularly useful in the study of synthetic polymers for shock absorption. Each model uses a distinct arrangement of these elements.
To analyze these features, a memory term must be introduced into the equation of  \eqref{eq:ondaI}. This is mathematically represented by a convolution in time, that is
\begin{equation}\label{eq:ondakI}
u_{tt} -\triangle u+ \int_0^t k(t-s) \triangle u(s)ds=0.
\end{equation}
A weak solution of \eqref{eq:ondakI} is  a function $u\in C([0,\infty);H^1_0(\Omega))\cap C^1([0,\infty);L^2(\Omega))$ such that for any $v\in H^1_0\Omega)$ $ t\to\int_\Omega u_t v\ dx\in C^1([0,\infty))$ 
and
\begin{equation*}
\frac {d}{dt} \int_\Omega u_t v\ dx -\int_\Omega \nabla u\cdot\nabla v\ dx
+\int_\Omega\int_0^t k(t-s)\nabla u(s)\,ds\cdot\nabla v\ dx=0\,,
\qquad\forall t\ge0\,.
\end{equation*}
Assuming that the integral kernel $k: [0,\infty)\to [0,\infty)$ is a locally absolutely continuous function such that $k(0)> 0$ and $k'(t)\le 0$ for a.e. $t\ge 0$ 
the system is dissipative. Indeed, 
the energy  of a weak solution $u$ is defined by the formula
\begin{equation*}\label{eqn:energy}
\begin{split}
E(t)=&\frac{1}{2}\int_{\Omega} |u_t(t,x)|^2\ dx+\frac{1}{2}\Big(1-\int_0^t\ k(s) \ ds\Big)\int_{\Omega} |\nabla u(t,x)|^2\ dx
\\
&+
\frac{1}{2}\int_{\Omega}\int_0^t k(t-s) |\nabla u(s,x)-\nabla u(t,x)|^2\ ds\ dx\qquad t\ge0\,,
\end{split}
\end{equation*}
and $E(t)$ is a decreasing function, since
\begin{equation*}
E'(t)=
\frac{1}{2}\int_0^t k'(t-s)\int_{\Omega} |\nabla u(s,x)-\nabla u(t,x)|^2\ dx\ ds
-\frac12k(t)\int_{\Omega} |\nabla u(t,x)|^2\ dx\leq 0
\,,
\end{equation*}
see e.g. \cite{MS}.

Due to the dissipation from the viscoelastic memory, combined with the exponential decay of the integral kernel and the condition 
\begin{equation}\label{eq:int-k}
\int_0^{\infty} k (t)\ dt<1\,,
\end{equation}
exponential stability is guaranteed, see e.g. \cite{Berrimi,ACS,CS3}.

Moreover, in \cite{LS1} we prove the  so-called  direct inequality or hidden regularity,
 that is  the regularity for the normal derivatives 
 of the weak solutions  of \eqref{eq:ondakI}.

The behavior of viscoelastic materials is characterized through specific kernel functions $k$.
For example, the Maxwell model accurately describes most polymers, where  
the kernel function is given by
\begin{equation*}
k(t)=b e^{-r t},
\quad 0<b<r.
\end{equation*}
The condition $b<r$ guarantees that \eqref{eq:int-k} holds.

The study of some viscoelastic materials can be conducted by generalizing the Maxwell model and using Prony sums defined by
\begin{equation*}
k(t)=\sum_{i=1}^N b_i e^{- r_i t}
\quad b_i\,, r_i>0,
\quad i=1,\dots N.
\end{equation*}
It is significant to observe that the Prony sums $k(t)$ are positive functions such that 
$k(0)> 0$ and $k'(t) \leq 0 $,
since
\begin{equation*}
k(0)=\sum_{i=1}^N b_i >0,
\quad
k'(t)=-\sum_{i=1}^N b_i  r_i e^{- r_it}<0.
\end{equation*}
Since
$
\int_0^{\infty} k (t)\ dt=\sum_{i=1}^N\frac {b_i}{r_i},
$
assuming
\begin{equation}\label{eq:cond<}
\sum_{i=1}^N\frac {b_i}{r_i}<1,
\end{equation}
condition \eqref{eq:int-k} holds. 

On the contrary, for  glassy type materials the condition \eqref{eq:int-k}, and hence \eqref{eq:cond<}, is not satisfied, because in that case one has $\sum_{i=1}^N\frac {b_i}{r_i}=1$, see \cite{MM} and also \cite{LS0}. Few results are known in the case $\int_0^{\infty} k (t)\ dt=1$.
In \cite{LS00} we consider the Burger model, which combines the Maxwell and Kelvin-Voigt models, a significant tool in the study of glass relaxation. It approximates the stretched exponential function, often used to describe this process, through the use of Prony sums with $N=2$.
Taking inspiration from the Burger model, we prove the direct inequality under a weaker condition on a general kernel $k$,  assuming that 
$\int_0^{T} k (t)\ dt<1$ for any $T>0$.

In this paper, we focus on glassy kernels that generalize the features of Prony sums, specifically those satisfying the condition 
$\sum_{i=1}^N\frac {b_i}{r_i}=1$. More precisely, we consider locally absolutely continuous
functions $k: [0,\infty)\to [0,\infty)$ such that 
\begin{equation*}
k(0)> 0,
\qquad
k'(t)\le 0\quad\mbox{for a.e.}\,\, t\ge 0,
\end{equation*}
\begin{equation*}
\int_0^{\infty} k (t)\ dt=1.
\end{equation*}
We want to demonstrate the energy decay of the weak solution and estimate the decay constant, assuming the kernel exhibits exponential decay. This is characterized by the condition
\begin{equation*}
k'(t)\le-\eta\ k(t)\quad\mbox{for a.e.}\,\, t\ge 0 \ \ (\eta>0).
\end{equation*}
Since we assume that $\int_0^{\infty} k (t)\ dt=1$, we must employ a different strategy than the one used in papers (see e.g. \cite{Berrimi,ACS,CS3}) where the condition is 
 $\int_0^{\infty} k (t)\ dt<1$. In this way, we also obtained effective estimates that allowed us to estimate the decay constant, see \eqref{eq:const-dec} below.

We prove the following result.
\begin{theorem}
There exists a positive number $\alpha$ 
such that
the energy $E(t)$ of a weak solution  of \eqref{eq:ondakI} decays exponentially as 
\begin{equation*}
E(t)\le E(0)e^{1-\alpha t}\qquad\forall t\ge 0\,.
\end{equation*}
\end{theorem}
The proof is given in an abstract setting, where $A$ is a self-adjoint linear
  operator on a Hilbert space $H$ with dense domain $D(A)$, satisfying
\begin{equation*}
\langle Ax,x\rangle\ge C\|x\|^2\qquad\forall x\in D(A)
\end{equation*}
for some $C>0$.

The plan of the paper is the following. In Section 2 we give some preliminaries useful in the sequel.
Section 3 is devoted to state and prove the decay result. 

\section{Preliminaries}

Let $H$ be a real Hilbert space with scalar product
$\langle \cdot \, ,\, \cdot \rangle$ and norm $\| \cdot \|$.

Throughout this paper, we will assume that the integral kernel $k$ satisfies  the following conditions:
\begin{enumerate}
\item[]
\begin{equation}\label{eq:k1}
k: [0,\infty)\to [0,\infty)\ \mbox{ is a locally absolutely continuous function}
\qquad
k(0)> 0,
\end{equation}
\item[]
\begin{equation}\label{eq:k2}
\int_0^{\infty} k (t)\ dt=1,
\end{equation}
\item[]
\begin{equation}\label{eq:k3}
k'(t)\le-\eta\ k(t)\quad\mbox{for a.e.}\,\, t\ge 0 \ \ (\eta>0)\,.
\end{equation}
\end{enumerate}
Thanks to \eqref{eq:k3}, we note that
\begin{equation}\label{eq:decayk}
\int_t^\infty k(s)\ ds\le -\frac1\eta\int_t^\infty k'(s)\ ds=\frac1\eta k(t)\,.
\end{equation}
We recast in an abstract formulation
\begin{equation}\label{eq:stato}
u''(t) +Au (t)-\int_0^t k(t-s)Au(s)\,ds=0\,,
\qquad t\in (0,\infty)\,,
\end{equation}
where $A$ is a self-adjoint linear
  operator on $H$ with dense domain $D(A)$, satisfying
\begin{equation}\label{eq:opA}
\langle Ax,x\rangle\ge C\|x\|^2\qquad\forall x\in D(A)
\end{equation}
for some $C>0$.

The energy  of a weak solution $u$ of \eqref{eq:stato} is 
defined by
\begin{equation*}
\begin{split}
E(t)=&\frac{1}{2}\| u'(t)\|^2+\frac{1}{2}\Big(1-\int_0^t\ k(s) \ ds\Big)\|\ra u(t)\|^2
\\
&+\frac{1}{2}\int_0^t k(t-s)\|\ra u(s)-\ra u(t)\|^2\ ds\qquad t\ge0\,.
\end{split}
\end{equation*}
Thanks to the assumption $\int_0^\infty k(s)\ ds=1$, we have  $1-\int_0^t k(s)\ ds=\int_t^\infty k(s)\ ds$, so the expression of the energy  can be given in an equivalent way
\begin{equation}\label{eq:energy}
\begin{split}
E(t)=&\frac{1}{2}\| u'(t)\|^2+\frac{1}{2}\int_t^\infty k(s)\ ds\ \|\ra u(t)\|^2
\\
&+\frac{1}{2}\int_0^t k(t-s)\|\ra u(s)-\ra u(t)\|^2\ ds\qquad t\ge0\,.
\end{split}
\end{equation}
$E(t)$ is a decreasing function, since we have
\begin{equation}\label{eq:deriv_energy}
E'(t)=
\frac{1}{2} \int_0^t k'(t-s)\|\ra u(s)-\ra u(t)\|^2\ ds
-\frac12k(t)\|\ra u(t)\|^2
\qquad t\ge0
\,.
\end{equation}
For reader's convenience we recall a well-known result, see e.g. \cite[Theorem 8.1]{K}.
\begin{theorem}\label{th:exp_th}
Let $E$ be a nonnegative decreasing function defined on $[0,\infty)$. If for  some constant $\alpha> 0$ we have
\begin{equation}\label{eq:exp_th}
\alpha\int_S^\infty\,E(t)\,dt\le E(S)\qquad\forall S\ge 0,
\end{equation}
then
$$
E(t)\le E(0)e^{1-\alpha t}\qquad\forall t\ge 0\,.
$$
\end{theorem}
Throughout the paper we will  use a standard notation for the integral convolution between two functions, that is 
\begin{equation*}
v*w(t)=\int_0^tv(t-s) w(s)\ ds
\,.
\end{equation*}

\section{Main result}
\begin{theorem}\label{th:decayen}
Assume \eqref{eq:k1}, \eqref{eq:k2}, \eqref{eq:k3} and \eqref{eq:opA}.
Then, 
the energy $E(t)$ of the weak solution $u$ of \eqref{eq:stato} decays as 
\begin{equation}\label{eq:decay}
E(t)\le E(0)e^{1-\alpha t}\qquad\forall t\ge 0\,,
\end{equation}
where the decay constant is given by
\begin{equation}\label{eq:const-dec}
\alpha=\frac1{2\big(\frac{k(0)+2}{ C}+1+\frac3\eta\big)}.
\end{equation}

\end{theorem}
The above theorem is an immediate consequence of Theorem \ref{th:exp_th} and of the following
technical result.
\begin{theorem}\label{th:intin}
Let assumptions \eqref{eq:k1}-\eqref{eq:k3} and \eqref{eq:opA} be satisfied. 
Then,  there exists a positive number 
$\alpha$ such that for any $u_0\in D(\ra)$ and 
$u_1\in X$
the energy $E(t)$ of the weak solution $u$ of \eqref{eq:stato} satisfies
\begin{equation}\label{eqn:Komornik}
\alpha\int_S^\infty  E(t)\ dt\le \ E(S)\qquad \forall\,  S\ge 0\,.
\end{equation}
\end{theorem}
\subsection{Proof of Theorem \ref{th:intin}}
First, we need to establish a useful identity.
\begin{lemma}\label{le:iden1}
For any $T\ge S\ge 0$ the following identity holds.
 \begin{multline}\label{eq:iden1}
\int_S^T\ \|u'(t) \|^2\ dt  
\\   
=
\int_S^T   \langle u'(t) ,(k*u)'(t)\rangle\ dt
+\int_S^T \|\ra u (t)-k *\ra u(t)\|^2\ dt
+ \Big[\langle u'(t) ,u(t)-k*u(t)\rangle\Big]_S^T
\,.  
\end{multline}
\end{lemma}
\begin{Proof}
Let us take the scalar product of both sides of equation \eqref{eq:stato} with $ u-k*u$ and integrate over $[S,T]$.  We obtain
\begin{equation*}
\int_S^T  \ \langle u''(t) +Au (t)-k *Au(t),u(t)-k*u(t)\rangle
\ dt=0\,.
\end{equation*}
Integrating by parts we get
\begin{equation*}
\Big[ \langle u'(t),u(t)-k * u(t)\rangle\Big]_S^T-\int_S^T  \|u'(t)\|^2\ dt+\int_S^T   \langle u'(t) ,(k*u)'(t)\rangle\ dt
+ \int_S^T \|\ra u (t)-k *\ra u(t)\|^2\ dt=0
\,,
\end{equation*}
that is \eqref{eq:iden1}.
\end{Proof}
\begin{Proofy}
To prove \eqref{eqn:Komornik}, we will show that there exist a positive number 
$\alpha$ such 
the energy $E(t)$ given by  \eqref{eq:energy} satisfies
\begin{equation}\label{eqn:Komornik1}
\alpha\int_S^T  E(t)\ dt\le \ E(S)\qquad \forall\,  T\ge S\ge 0\,.
\end{equation}
To this end, we note that
\begin{equation}\label{eq:intenergy}
\begin{split}
\int_S^T E(t)\ dt=&\frac{1}{2}\int_S^T \| u'(t)\|^2\ dt+\frac{1}{2}\int_S^T \int_t^\infty k(s)\ ds\ \|\ra u(t)\|^2\ dt
\\
&+\frac{1}{2}\int_S^T \int_0^t k(t-s)\|\ra u(s)-\ra u(t)\|^2\ ds\ dt\qquad t\ge0\,.
\end{split}
\end{equation}
First, we estimate the second and the third term on the right-hand side of \eqref{eq:intenergy}. Indeed, by \eqref{eq:decayk} we get
\begin{equation*}
 \int_t^\infty k(s)\ ds\ \|\ra u (t)\|^2
\le\frac1\eta k(t)\|\ra u(t)\|^2,
\end{equation*}
and by 
\begin{equation*}
\int_0^tk(t-s)\|\ra u(s)-\ra u(t)\|^2\ ds
\le
-\frac{1}{\eta}\int_0^tk'(t-s)\|\ra u(s)-\ra u(t)\|^2\ ds
\,,
\end{equation*}
and hence by \eqref{eq:deriv_energy} 
\begin{equation}\label{eq:2-3energy}
\frac{1}{2} \int_t^\infty k(s)\ ds\ \|\ra u (t)\|^2+ \frac{1}{2} \int_0^tk(t-s)\|\ra u(s)-\ra u(t)\|^2\ ds
\le
-\frac{1}{\eta}E'(t)\,.
\end{equation}
Therefore
\begin{multline}\label{eq:3energy} 
\frac{1}{2}\int_S^T \Big(\int_t^\infty k(s)\ ds\ \|\ra u (t)\|^2+  \int_0^tk(t-s)\|\ra u(s)-\ra u(t)\|^2\ ds\Big)\ dt
\\
\le
-\frac{1}{\eta}\int_S^TE'(t)\ dt\le \frac{1}{\eta} E(S)\,.
\end{multline}
To evaluate the term $\int_S^T \| u'(t)\|^2\ dt$ we use Lemma \ref{le:iden1}: we have to estimate the terms on the right-hand side of \eqref{eq:iden1}. Beginning with the first term, we note that
\begin{equation}\label{eq:epsilon}
\int_S^T   \langle u'(t) ,(k*u)'(t)\rangle\ dt
\le\frac{\varepsilon}2\int_S^T \|u'(t)\|^2\ dt+\frac1{2\varepsilon}\int_S^T \|(k*u)'(t)\|^2\ dt
\end{equation}

\begin{equation*}
\begin{split}
(k*u)'(t)
=&\int_0^t\ k'(t-s) u(s)\ ds+k(0)u(t)
\\
=&\int_0^t\ k'(t-s) (u(s)-u(t))\ ds+k(t)u(t)
\,,
\end{split}
\end{equation*}
so we have
\begin{equation*}
\frac12\|(k*u)'(t)\|^2
\le
\Big\|\int_0^t\ k'(t-s) (u(s)- u(t))ds\Big\|^2
+k(t)^2\|u(t)\|^2.
\end{equation*}
Since
\begin{equation*}
\begin{split}
&\Big\|\int_0^t\ k'(t-s) (u(s)- u(t))ds\Big\|^2  
\\
\le&
\Big(\int_0^t\ |k'(t-s)|^{1/2} |k'(t-s)|^{1/2}\| u(s)- u(t)\|ds\Big)^2  
\\
\le&
\int_0^t\ |k'(s)|\ ds\int_0^t\  |k'(t-s)|\| u(s)- u(t)\|^2 ds
\\
=&
\big(k(t)-k(0)\big)\int_0^t\  k'(t-s)\| u(s)- u(t)\|^2 ds
\\
\le&
-k(0)\int_0^t\  k'(t-s)\| u(s)- u(t)\|^2 ds
\,,
\end{split}
\end{equation*}
and
\begin{equation*}
k(t)^2\|u(t)\|^2
\le 
k(0)k(t)\|u(t)\|^2,
\end{equation*}
taking into account \eqref{eq:opA} and \eqref{eq:deriv_energy} we get
\begin{equation*}
\begin{split}
&\frac12\|(k*u)'(t)\|^2
\\
\le&
\frac{2k(0)}{C}\Big(-\frac12\int_0^t\  k'(t-s)\|\ra u(s)- \ra u(t)\|^2 ds
+\frac12k(t)\|\ra u(t)\|^2\Big)
\\
=&-\frac{2k(0)}{C}E'(t),
\end{split}
\end{equation*}
and hence
\begin{equation*}
\frac12\int_S^T\|(k*u)'(t)\|^2\ dt\le -\frac{2k(0)}{C}\int_S^TE'(t)\ dt\le \frac{2k(0)}{ C}E(S)\,.
\end{equation*}
Plugging the above estimate into \eqref{eq:epsilon} we obtain
\begin{equation}\label{eq:firstterm}
\int_S^T   \langle u'(t) ,(k*u)'(t)\rangle\ dt
\le\frac{\varepsilon}2\int_S^T \|u'(t)\|^2\ dt+\frac{2k(0)}{\varepsilon C}E(S)
\le\varepsilon\int_S^T E(t)\ dt+\frac{2k(0)}{\varepsilon C}E(S)\,.
\end{equation}
Now, since $\int_0^\infty k(s)\ ds=1$, we have
\begin{multline}\label{eq:u-k*u}
u(t)-k*u(t)
=\int_t^\infty k (s)\ ds\ u(t)+\int_0^t k (s)\ ds\ u(t)-\int_0^t k (t-s)u(s) ds
\\
=\int_t^\infty k (s)\ ds\ u(t)-\int_0^t k (t-s)(u(s)-u(t)) ds\,,
\end{multline}
and hence
\begin{equation*}
\| \ra u-k *\ra u\|^2
\le
2\Big(\int_t^\infty k (s) ds\Big)^2 \|\ra u(t)\|^2
+2\Big(\int_0^t k (t-s)\|\ra u(s)-\ra u(t)\|ds\Big)^2\,.
\end{equation*}
We observe that
\begin{equation*}
\int_0^t k(t-s)\  \|\ra u(s)-\ra u(t) \| ds 
  \le     
\Big(\int_0^{t} k(s) ds\Big)^{1/2}\Big(\int_0^tk(t-s)\|\ra u(s)-\ra u(t))\|^2\ ds\Big)^{1/2}
\,,
\end{equation*}
so, again by assumption $\int_0^\infty k(s)\ ds=1$, we get
\begin{equation}\label{eq:energy-1}
\| \ra u-k *\ra u\|^2
\le
2\int_t^\infty k (s) ds\  \|\ra u(t)\|^2
+2\int_0^t k (t-s)\|\ra u(s)-\ra u(t)\|^2ds\,.
\end{equation}
Using \eqref{eq:2-3energy}
we have
\begin{equation*}
\| \ra u-k *\ra u\|^2
\le
-\frac{4}{\eta}E'(t)\,,
\end{equation*}
whence
\begin{equation}\label{eq:secondterm}
\int_S^T\| \ra u-k *\ra u\|^2\ dt
\le
-\frac{4}{\eta}\int_S^TE'(t)\ dt\le\frac{4}{\eta}E(S)\,,
\end{equation}
so the estimate of the second term on the right-hand side of \eqref{eq:iden1} is completed.

Finally, we have to estimate the last term $\Big[\langle u'(t) ,u(t)-k*u(t)\rangle\Big]_S^T$ on the right-hand side of \eqref{eq:iden1}.
We note that 
\begin{equation*}
|\langle u'(t),u(t)-k*u(t)\rangle|
\le \frac{1}{2}\|u'(t)\|^2+\frac{1}{2}\|u(t)-k*u(t)\|^2\,.
\end{equation*}
By \eqref{eq:opA} and  \eqref{eq:energy-1} we obtain
\begin{multline*}
\frac{1}{2}\|u(t)-k*u(t)\|^2
\le
\frac{1}{2C}\|\ra u(t)-k*\ra u(t)\|^2
\\ 
\le
\frac{1}{C} \int_t^\infty k (s) ds\  \|\ra u(t)\|^2
+\frac{1}{C}\int_0^t k (t-s)\|\ra u(s)-\ra u(t)\|^2ds\,,
\end{multline*}
and hence by \eqref{eq:energy}
\begin{equation*}
|\langle u'(t),u(t)-k*u(t)\rangle|
\le \Big(1+\frac{2}{C}\Big) E(t)\,.
\end{equation*}
Since $E(t)$ is decreasing, we have
\begin{equation}\label{eq:thirdterm}
\Big[\langle u'(t) ,u(t)-k*u(t)\rangle\Big]_S^T
\le2\Big(1+\frac{2}{C}\Big)E(S)\,.
\end{equation}
In conclusion, by \eqref{eq:firstterm}, \eqref{eq:secondterm}, \eqref{eq:thirdterm} and \eqref{eq:3energy} we get
\begin{equation*}
\int_S^T E(t)\ dt
\le \frac\varepsilon2\int_S^T E(t)\ dt+\Big(\frac{k(0)}{\varepsilon C}+1+\frac{2}{C}+\frac3\eta\Big)E(S)\,.
\end{equation*}
For $\varepsilon=1$ we have
\begin{equation*}
\int_S^T E(t)\ dt
\le 2\Big(\frac{k(0)}{ C}+1+\frac{2}{C}+\frac3\eta\Big)E(S)\,,
\end{equation*}
so, taking $\alpha=\frac1{2\big(\frac{k(0)+2}{ C}+1+\frac3\eta\big)}$, the proof of \eqref{eqn:Komornik1} has been completed.
\end{Proofy}

\begin{remark}
If $\Omega\subset\R^N$ $(N\ge1)$ is a bounded open set, we define the operator $A$ in $L^2(\Omega)$ by
\begin{equation*}
\begin{split}
D(A)&= H^2(\Omega)\cap H^1_0(\Omega)
\\
(Au)(x)&=-\triangle u(x),
\quad u\in D(A), \  x\in\Omega .
\end{split}
\end{equation*}
In this specific case
we can express the decay constant \eqref{eq:const-dec} in a more explicit form. 
Since the spectrum of $A$ consists of a sequence of positive eigenvalues tending to $+\infty$, denoting  the lowest nonzero eigenvalue of $A$ with $\lambda_1$, the constant $C$ assumed in \eqref{eq:opA} simplifies to 
$C=\lambda_1$, and hence $\alpha=\frac1{2\big(\frac{k(0)+2}{ \lambda_1}+1+\frac3\eta\big)}$.
\end{remark}

\end{document}